\documentclass[a4paper,10pt]{amsart}

\usepackage[cp850]{inputenc}
\usepackage{marvosym}
\usepackage{wasysym}
\usepackage{latexsym}
\usepackage{amsfonts}
\usepackage{amsmath}
\usepackage{amssymb,amsxtra,amscd}
\usepackage{mathrsfs}
\usepackage{hyperref}

\newcommand{\R}{\mathbb{R}}
\newcommand{\C}{\mathbb{C}}
\newcommand{\N}{\mathbb{N}}

\newcommand{\supp}{\mbox{supp}\,}
\newcommand{\singsupp}{\mbox{sing\,supp}\,}
\newcommand{\dist}{\mbox{dist}}

\newtheorem{theo}{Theorem}
\newtheorem{lem}[theo]{Lemma}

\newtheorem{prop}[theo]{Proposition}

\newtheorem{rem}[theo]{Remark}

\newtheorem{example}[theo]{Example}

\title{Every $P$-convex subset of $\R^2$ is already strongly $P$-convex}
\author{T. Kalmes}

\begin{document}

\maketitle

\begin{center}
{\sl Dedicated to the memory of Susanne Dierolf}
\end{center}

\begin{abstract}
A classical result of Malgrange says that for a polynomial $P$ and an open subset $\Omega$ of $\R^d$ the differential operator $P(D)$ is surjective on $C^\infty(\Omega)$ if and only if $\Omega$ is $P$-convex. H\"ormander showed that $P(D)$ is surjective as an operator on $\mathscr{D}'(\Omega)$ if and only if $\Omega$ is strongly $P$-convex. It is well known that the natural question whether these two notions coincide has to be answered in the negative in general. However, Tr\`eves conjectured that in the case of $d=2$ $P$-convexity and strong $P$-convexity are equivalent. A proof of this conjecture is given in this note.
\end{abstract}

\section{Introduction}

It is a classical result by Malgrange \cite[Chapitre 1, Th\'eor\`eme 4]{Malgrange} that for a polynomial $P\in\C[X_1,\ldots,X_d]$ and for an open set $\Omega\subset\R^d$ the constant coefficient differential operator $P(D):C^\infty(\Omega)\rightarrow C^\infty(\Omega)$ is surjective if and only if $\Omega$ is $P$-convex, that is, if and only if for every compact subset $K$ of $\Omega$ there is another compact subset $L$ of $\Omega$ such that for each $u\in\mathscr{E}'(\Omega)$ with $\supp P(-D)u\subset K$ it holds $\supp u\subset L$.

H\"ormander showed \cite{Hoermander 1} that $P(D)$ is surjective as an operator on $\mathscr{D}'(\Omega)$ if and only if $\Omega$ is strongly $P$-convex, i.e.\ $\Omega$ is $P$-convex as well as $P$-convex for singular supports, the later meaning that for every compact subset $K$ of $\Omega$ there is another compact subset $L$ of $\Omega$ such that for each $u\in\mathscr{E}'(\Omega)$ with $\singsupp P(-D)u\subset K$ it holds $\singsupp u\subset L$.

Clearly, stong $P$-convexity implies $P$-convexity and it is a natural question to ask if (or when) these notions coincide. It is well-known that in general the answer to this question is in the negative. However, Tr\`eves conjectured \cite[p.\ 389, Problem 2]{Treves} that in the case of $\Omega\subset\R^2$, $P$-convexity and strong $P$-convexity are equivalent, i.e.\ for an open subset $\Omega$ of $\R^2$ surjectivity of $P(D):C^\infty(\Omega)\rightarrow C^\infty(\Omega)$ is equivalent to surjectivity of $P(D):\mathscr{D}'(\Omega)\rightarrow\mathscr{D}'(\Omega)$.\\

From now on we will use the terminology of \cite{Hoermander}. In particular, we call $P$-convexity for supports what is called $P$-convexity above. Hence we will have proved Tr\`eves conjecture if we prove the following theorem.

\begin{theo}\label{treves conjecture}
Let $\Omega\subset\R^2$ be open and $P\in\C[X_1,X_2]$. If $\Omega$ is $P$-convex for supports then $\Omega$ is already $P$-convex for singular supports.
\end{theo}

In order to prove Theorem \ref{treves conjecture} we will apply H\"ormander's theory of continuation of differentiability (cf.\ \cite[Section 11.3., vol.\ II]{Hoermander}).\\

The paper is organized as follows. In section 2 we will expose the connection of the localizations at infinity of a polynomial $P$ and a certain real-valued function $\sigma_P$ defined on the subspaces of $\R^d$. This will help us to see that in case of $d=2$ for a given $P$ certain important hyperplanes are always characteristic. In section 3 we will give sufficient conditions on an open subset $\Omega$ of $\R^d$ to be $P$-convex for supports as well as $P$-convex for singular supports. These will be applied in section 4 in order to prove Theorem \ref{treves conjecture}.

Throughout the paper we use standard notation from distribution theory and partial differential operators as may be found in \cite{Hoermander}. In order to avoid cumbersome formulations we assume that $P$ is non-zero throughout the whole paper. Moreover, for a hyperplane $H=\{x\in\R^d;\,\langle x,N\rangle=\alpha\}$ with $N\in S^{d-1},\alpha\in\R,$ we denote by $H^\perp$ the linear span of $N$.

\section{Localizations at Infinity and Continuation of Differentiability}

The problem we want to solve is clearly related to deriving bounds for $\singsupp u$ by knowledge of $\singsupp P(-D)u$, where $u\in\mathscr{E}'(\Omega)$ for $\Omega\subset\R^d$ open. If $E$ is a fundamental solution of $\check{P}$ we have $u=P(-D)u* E$ and from this it follows that for the wave front set $WF(u)$ of $u$ one has
\begin{eqnarray}
WF(u)\subset\{(x+y,\xi);\, (x,\xi)\in WF(P(-D)u)\mbox{ and }(y,\xi)\in WF(E)\}
\end{eqnarray}
(cf.\ \cite[p.\ 270, vol.\ I, Formula (8.2.16)]{Hoermander}), where the wave front set of a distribution $v$ is a subset of $\R^d\times S^{d-1}$ whose projection onto $\R^d$ is precisely $\singsupp v$. Therefore, knowledge about $WF(P(-D)u)$ as well as $WF(E)$ will allow to obtain bounds for $\singsupp u$.

For every polynomial $P$ there is a specific fundamental solution $E(P)$ for which the location of its wave front set is well understood by means of the so called localizations at infinity of $P$ whose definition we want to recall.

For a polynomial $P$ and $\xi\in\R^d$ we set $P_\xi(\eta)=P(\eta+\xi)$. The set of limits of the normalized polynomials
\[\eta\mapsto\frac{P_\xi(\eta)}{\tilde{P_\xi}(0)}\]
as $\xi$ tends to infinity is denoted by $L(P)$, where $\tilde{P_\xi}(0)=\sqrt{\sum_\alpha |P^{(\alpha)}_\xi(0)|^2}$ and where for a multiindex $\alpha\in\N_0^d$ we denote the $\alpha$-derivative of $P_\xi$ by $P^{(\alpha)}_\xi$. More precisely, if $N\in S^{d-1}$ then the set of limits where $\xi/|\xi|\rightarrow N$ is denoted by $L_N(P)$. Obviously, $L(P)$ as well as $L_N(P)$ are closed subsets of the unit sphere of all polynomials in $d$ variables of degree not exceeding the degree of $P$, equipped with the norm $Q\mapsto\tilde{Q}(0)$. The non-zero multiples of elements of $L(P)$ (resp.\ of $L_N(P)$) are called {\sl localizations of $P$ at infinity} (resp.\ {\sl localizations of $P$ at infinity in direction $N$}). Clearly, $Q\in L_{N}(\check{P})$ if and only if $\check{Q}\in L_{-N}(P)$.

Recall that for a polynomial $Q$
\[\Lambda(Q)=\{\eta\in\R^d;\forall\xi\in\R^d,t\in\R:\,Q(\xi+t\eta)= Q(\xi)\},\]
which is obviously a subspace of $\R^d$. Moreover, denote by $\Lambda'(Q)$ the orthogonal space of $\Lambda(Q)$. Clearly, $Q$ is constant if and only if $\Lambda'(Q)=\{0\}$. By a result due to H\"ormander (cf.\ \cite[Theorem 10.2.11, vol.\ II]{Hoermander}) the wave front set $WF(E(\check{P}))$ of the above mentioned fundamental solution $E(\check{P})$ is contained in the closure of the set
\[\{(x,N)\in\R^d\times S^{d-1};\ x\in\Lambda'(Q)\mbox{ for some }Q\in L_N(\check{P})\}.\]
From this it clearly follows that for $u\in\mathscr{E}'(\Omega)$ the non-constant elements of $L(\check{P})$ are the ones which may cause $\singsupp u$ to be much larger than $\singsupp P(-D)u$ due to equation (1) above.

Define for a polynomial $Q$, a subspace $V$ of $\R^d$, and $t\geq 1$
\[\tilde{Q}_V(\xi,t)=\sup\{|Q(\xi+\eta)|;\,\eta\in V,|\eta|\leq t\}\]
and
\[\tilde{Q}(\xi,t)=\tilde{Q}_{\R^d}(\xi,t).\]
Clearly, for every $\xi\in\R^d$ and $t\geq 1$ $\tilde{Q}(\xi,t)$ is a norm on the space of all polynomials.
So, if $Q\in L(\check{P})$ is non-constant then
\[0=\inf_{t\geq 1}\frac{\tilde{Q}_{\Lambda(Q)}(0,t)}{\tilde{Q}(0,t)}\]
because the numerator equals $|Q(0)|$ while the denominator tends to infinity with $t$. Moreover, since $Q\in L(\check{P})$ it follows that there is a sequence $(\xi_n)_{n\in\N}$ in $\R^d$ tending to infinity such that $Q=\lim_{n\rightarrow}\check{P}_{\xi_n}/\tilde{\check{P}}_{\xi_n}(0)$, hence
\[0=\inf_{t\geq 1}\frac{\tilde{Q}_{\Lambda(Q)}(0,t)}{\tilde{Q}(0,t)}=\inf_{t\geq 1}\lim_{n\rightarrow\infty}\frac{\tilde{\check{P}}_{\Lambda(Q)}(\xi_n,t)}{\tilde{\check{P}}(\xi_n,t)}.\]

Defining for an arbitrary subspace $V$ of $\R^d$
\[\sigma_{\check{P}}(V)=\inf_{t\geq 1}\liminf_{\xi\rightarrow\infty}\frac{\tilde{\check{P}}_V(\xi,t)}{\tilde{\check{P}}(\xi,t)},\]
it follows immediately that $\sigma_{\check{P}}(V)=\sigma_P(V)$. Moreover, for $y\in\R^d$ we shall simply write $\sigma_P(y)$ instead of $\sigma_P(span\{y\})$. The function $\sigma_P$ is much more powerful than simply identifying non-constant elements of $L(\check{P})$.

The values of $\sigma_P$ govern the possibility to continue differentiability of zero solutions of $P(D)$ across a hyperplane $H=\{x;\langle x,N\rangle=\alpha\}, N\in S^{d-1},\alpha\in\R$: Let $\Omega\subset\R^d$ be open, $x_0\in\Omega$ and $N\in S^{d-1}$ be such that $\sigma_P(N)\neq 0$. Then there is a neighborhood $U$ of $x_0$ such that $u\in C^\infty(U)$ for every $u\in\mathscr{D}'(\Omega)$ with $P(D)u=0$ as well as $u_{|\Omega_-}\in C^\infty(\Omega_-)$, where $\Omega_-=\{x\in\Omega;\,\langle x,N\rangle<\langle x_0,N\rangle\}$. This is only a very special case of \cite[Theorem 11.3.6, vol.\ II]{Hoermander}.

We have already indicated the connection between the localizations of $P$ at infinity and the function $\sigma_P$. The next lemma contains some more results which will be needed in the sequel.

\begin{lem}\label{sigmaP via localizations}
Let $P$ be of degree $m$ with principal part $P_m$.
\begin{itemize}
    \item[i)] For every subspace $V$ of $\R^d$ and $t\geq 1$ we have \[\liminf_{\xi\rightarrow\infty}\frac{\tilde{P}_V(\xi,t)}{\tilde{P}(\xi,t)}=\inf_{Q\in L(P)}\frac{\tilde{Q}_V(0,t)}{\tilde{Q}(0,t)}.\]
    \item[ii)] Let $N\in S^{d-1}$ and $Q\in L_N(P)$. If $P_m(N)\neq 0$ then $Q$ is constant.
    \item[iii)] If $P$ is non-elliptic then for every subspace $V$ of $\R^d$ and $t\geq 1$ we have
    \[\liminf_{\xi\rightarrow\infty}\frac{\tilde{P}_V(\xi,t)}{\tilde{P}(\xi,t)}=\inf_{N\in S^{d-1}, P_m(N)=0}\;\inf_{Q\in L_N(P)}\frac{\tilde{Q}_V(0,t)}{\tilde{Q}(0,t)}.\]
\end{itemize}
\end{lem}

{\sc Proof.} i) Since for every subspace $V$ and each $t\geq 1$ the maps $R\mapsto\tilde{R}_V(0,t)$ are continuous seminorms on the space of all polynomials $R$ in $d$ variables and because $\tilde{P}_V(\xi,t)=(\tilde{P_\xi})_V(0,t)$ it follows immediately from the definition that
\[\frac{\tilde{Q}_V(0,t)}{\tilde{Q}(0,t)}\geq \liminf_{\xi\rightarrow\infty}\frac{\tilde{P}_V(\xi,t)}{\tilde{P}(\xi,t)}\]
for every $Q\in L(P)$.

Moreover, if $(\xi_n)_{n\in\N}$ tends to infinity such that
\[ \liminf_{\xi\rightarrow\infty}\frac{\tilde{P}_V(\xi,t)}{\tilde{P}(\xi,t)}=\lim_{n\rightarrow\infty}\frac{\tilde{P}_V(\xi_n,t)}{\tilde{P}(\xi_n,t)}=\lim_{n\rightarrow\infty}\frac{(\tilde{P}_{\xi_n})_V(0,t)}{\tilde{P}_{\xi_n}(0,t)}\]
we can extract a subsequence of $(\xi_n)_{n\in\N}$ which we again denote by $(\xi_n)_{n\in\N}$ such that the sequence of normalized polynomials $P_{\xi_n}/\tilde{P}_{\xi_n}(0)$ converges in the compact unit sphere of all polynomials in $d$ variables of degree at most $m$. This limit belongs to $L(P)$ and we get
\[\liminf_{\xi\rightarrow\infty}\frac{\tilde{P}_V(\xi,t)}{\tilde{P}(\xi,t)}\geq \inf_{Q\in L(P)}\frac{\tilde{Q}_V(0,t)}{\tilde{Q}(0,t)}\]
completing the proof of i).

The proof of ii) is an easy application of Taylor's formula. Let $P=\sum_{j=0}^m P_j$, where $P_j$ is either a homogeneous polynomial of degree $j$ or identically zero. Let $(\xi_n)_{n\in\N}$ tend to infinity with $\lim_{n\rightarrow\infty}\xi_n/|\xi_n|=N$ and $P_m(N)\neq 0$. Then
\begin{eqnarray*}
    P_{\xi_n}(\eta)&=&\sum_{0\leq |\alpha|\leq j\leq m}\frac{P_j^{(\alpha)}(\xi_n)}{\alpha !}\eta^\alpha\\
    &=&|\xi_n|^m\left(\sum_{0\leq j\leq m}\frac{|\xi_n|^j}{|\xi_n|^m} P_j(\frac{\xi_n}{|\xi_n|})+\sum_{0 <|\alpha|\leq j\leq m}\frac{|\xi_n|^{j-|\alpha|}}{|\xi_n|^m\alpha !}P_j^{(\alpha)}(\frac{\xi_n}{|\xi_n|})\eta^\alpha\right).
\end{eqnarray*}
Moreover
\begin{eqnarray*}
    \tilde{P}_{\xi_n}(0)&=&\sqrt{\sum_{0\leq|\alpha|\leq m}|\sum_{j=|\alpha|}^mP_j^{(\alpha)}(\xi_n)|^2}\\
    &=&|\xi_n|^m\sqrt{|\sum_{j=0}^m P_j(\frac{\xi_n}{|\xi_n|})\frac{|\xi_n|^j}{|\xi_n|^m}|^2+\sum_{0<|\alpha|\leq m}|\sum_{j=|\alpha|}^m P_j^{(\alpha)}(\frac{\xi_n}{|\xi_n|})\frac{|\xi_n|^{j-|\alpha|}}{|\xi_n|^m}|^2},
\end{eqnarray*}
which implies that
\[\lim_{n\rightarrow\infty}\frac{P_{\xi_n}(\eta)}{\tilde{P}_{\xi_n}(0)}=\frac{P_m(N)}{|P_m(N)|}\]
for every $\eta\in\R^d$ showing ii).

iii) is an immediate consequence of i), ii), and $\liminf_{\xi\rightarrow\infty}\tilde{P}_V(\xi,t)/\tilde{P}(\xi,t)\leq 1$.\hfill$\square$

\begin{rem}\label{hypoelliptic remark}
\begin{rm}
Since for every localization $Q$ of $P$ at infinity one has $\Lambda(Q)\neq 0$ (cf.\ \cite[Theorem 10.2.8,vol.\ II]{Hoermander}) it follows that in case of $Q$ being non-constant there is a subspace $V\neq 0$ such that $\sigma_P(V)=0$. Recall that a polynomial $P$ is called {\it hypoelliptic} if $\singsupp P(D)u=\singsupp u$ for every $u\in\mathscr{D}'(\Omega)$, where $\Omega\subset\R^d$ is an arbitrary open set. As shown in the proof of \cite[Theorem 11.1.11, vol.\ II]{Hoermander} $P$ being hypoelliptic is equivalent to the fact that every localization of $P$ at infinity is constant. By the above lemma and the obvious fact that $\sigma_P(V_1)\leq\sigma_P(V_2)$ whenever $V_1\subset V_2$ it therefore follows easily that $P$ is hypoelliptic if and only if $\sigma_P(y)\neq 0$ for every $y\in\R^d$. Moreover, it is well-known that elliptic polynomials are hypoelliptic (cf.\ \cite[Theorem 11.1.10, vol.\ II]{Hoermander}).
\end{rm}
\end{rem}

The reason, why the case $d=2$ is so very different from the higher dimensional cases is because a non-zero homogeneous polynomial in two variables can only have a finite number of zeros in the unit sphere. With this observation we can prove the following key lemma.

\begin{lem}\label{singulars are characteristic in r2}
Let $P\in\C[X_1,X_2]$ be of degree $m$ with principal part $P_m$. Then
\[\{y\in S^1;\,\sigma_P(y)=0\}\subset\{y\in S^1;\,P_m(y)=0\}.\]
\end{lem}

{\sc Proof.} By Remark \ref{hypoelliptic remark} we can assume without loss of generality that $P$ is not hypoelliptic, hence not elliptic. Let $\{N\in S^1;\,P_m(N)=0\}=\{N_1,\ldots,N_l\}$. For each $1\leq j\leq l$ choose $x_j\in S^1$ orthogonal to $N_j$. Take an arbitrary, non-constant $Q\in L(P)$. By Lemma \ref{sigmaP via localizations} ii) there is $1\leq j\leq l$ such that $Q\in L_{N_j}(P)$. By \cite[Theorem 10.2.8, vol.\ II]{Hoermander} we have $Q(\xi+sN_j)=Q(\xi)$ for any $\xi\in\R^2, s\in\R$. Hence $Q(\xi)=Q(\langle\xi,x_j\rangle x_j)$ for all $\xi\in\R^2$. Defining
\[q:\R\rightarrow\C, s\mapsto Q(s x_j)\]
it follows that for fixed $y\in S^1$
\begin{eqnarray*}
    \tilde{Q}_{span\{y\}}(0,t)&=&\sup\{|Q(\lambda y)|;\,|\lambda|\leq t\}
    =\sup\{|Q(\lambda\langle y,x_j\rangle x_j)|;\,|\lambda|\leq t\}\\
    &=&\sup\{|q(\lambda t\langle y,x_j\rangle)|;\,|\lambda|\leq 1\},
\end{eqnarray*}
and because $|x_j|=1$ we also have
\begin{eqnarray*}
    \tilde{Q}(0,t)&=&\sup\{|Q(\xi)|;\,\xi\in\R^2, |\xi|\leq t\}
    =\sup\{|Q(\langle \xi,x_j\rangle x_j)|;\,\xi\in\R^2, |\xi|\leq t\}\\
    &=&\sup\{|Q(\lambda x_j)|;\,|\lambda|\leq t\}
    =\sup\{|q(\lambda t)|;\,|\lambda|\leq 1\}.
\end{eqnarray*}
Since $Q\in L(P)$ it follows that $q$ is a polynomial of degree at most $m$. Since on the finite dimensional space of all polynomials in one variable of degree at most $m$ the norms $\sup_{|s|\leq 1}|p(s)|$ and $\sum_{k=0}^m|p^{(k)}(0)|$ are equivalent there is $C>0$ such that
\[C\sup_{|s|\leq 1}|p(s)|\geq \sum_{k=0}^m|p^{(k)}(0)|\geq 1/C\,\sup_{|s|\leq 1}|p(s)|\]
for all $p\in\C[X]$ with degree at most $m$. Applying this to the polynomials $s\mapsto q(s t)$ and $s\mapsto q(s t\langle y,x_j\rangle)$ gives
\begin{eqnarray*}
    \frac{\tilde{Q}_{span\{y\}}(0,t)}{\tilde{Q}(0,t)}&\geq &\frac{\sum_{k=0}^m |q^{(k)}(0)|t^k|\langle y,x_j\rangle|^k}{C^2 \sum_{k=0}^m |q^{(k)}(0)|t^k}\\
    &\geq & |\langle y,x_j\rangle|^m/C^2,
\end{eqnarray*}
where we used $|\langle y,x_j\rangle|\leq 1$ in the last inequality. We conclude that for every $1\leq j\leq l$
\[\inf_{Q\in L_{N_j}(P)}\frac{\tilde{Q}_{span\{y\}}(0,t)}{\tilde{Q}(0,t)}\geq \frac{|\langle y,x_j\rangle|^m}{C^2},\]
where $C$ only depends on the degree $m$ of $P$. It follows from Lemma \ref{sigmaP via localizations} iii) and $\{N\in S^1;\, P_m(N)=0\}=\{N_1,\ldots, N_l\}$ that for all $t\geq 1$
\[\liminf_{\xi\rightarrow\infty}\frac{\tilde{P}_{span\{y\}}(\xi,t)}{\tilde{P}(\xi,t)}=\min_{1\leq j\leq l}\inf_{Q\in L_{N_j}(P)}\frac{\tilde{Q}_{span\{y\}}(0,t)}{\tilde{Q}(0,t)}\geq \min_{1\leq j\leq l}\frac{|\langle y,x_j\rangle|^m}{C^2}.\]
Therefore, if for $y\in S^1$
\[0=\sigma_P(y)=\inf_{t\geq 1}\liminf_{\xi\rightarrow\infty}\frac{\tilde{P}_{span\{y\}}(\xi,t)}{\tilde{P}(\xi,t)}\]
it follows that $y$ is orthogonal to some $x_j$, hence $y\in\{N_j,-N_j\}$ since $|y|=1=|N_j|$ which shows $P_m(y)=0$.\hfill$\square$\\

In particular we conclude that for $P\in \C[X_1,X_2]\backslash\{0\}$ the set $\{y\in S^1;\sigma_P(y)=0\}$ is finite. The next example shows that an analogous statement of the above lemma is not true in general in case of $d>2$.

\begin{example}\label{wave in rd}
\begin{rm}
Let $d>2$ and $P\in\C[X_1,\ldots,X_d]$ be given by
\[P(x_1,\ldots,x_d)=x_1^2-x_2^2-\ldots -x_d^2.\]
It follows that a localization of $P$ at infinity in direction $1/\sqrt{2}\,(1,1,0,\ldots,0)$ is given by $Q(\xi_1,\ldots,\xi_d)=(\xi_1-\xi_2)/2$. Hence it follows for $e_d=(0,\ldots,0,1)$ that $\tilde{Q}_{span\{e_d\}}(0,t)=0$ for every $t\geq 1$ so that $\sigma_P(e_d)=0$ by Lemma \ref{sigmaP via localizations}. On the other hand, we clearly have $P_2(e_d)=P(e_d)=-1$.
\end{rm}
\end{example}

One way we will make use of $\sigma_P(V)$ is given by the following result which is nothing but a reformulation of \cite[Corollary 11.3.7, vol.\ II]{Hoermander}. For the proof see \cite[Corollary 3]{Frerick}.

\begin{prop}\label{Hoermander reformulated}
Let $\Omega_1\subset \Omega_2$ be open and convex subsets of $\R^d$, and let $P$ be a polynomial. Then the following are equivalent:
\begin{itemize}
    \item[i)] Every $u\in\mathscr{D}'(\Omega_2)$ satisfying $P(D)u\in C^\infty (\Omega_2)$ as well as $u|_{\Omega_1}\in C^\infty(\Omega_1)$ belongs to $C^\infty(\Omega_2)$.
    \item[ii)] Every hyperplane $H=\{x;\langle x,N \rangle=\alpha\}$ with $\sigma_P(N)=0$ which intersects $\Omega_2$ already intersects $\Omega_1$.
\end{itemize}
\end{prop}

It follows immediately from Lemma \ref{singulars are characteristic in r2} that in case of $d=2$ every hyperplane $H$ with $\sigma_P(H^\perp)=0$ is characteristic for $P$.

\section{Exterior Cone Conditions for $P$-convexity}

In this section we will prove some sufficient conditions for an open subset $\Omega$ of $\R^d$ to be $P$-convex for supports as well as $P$-convex for singular supports in terms of exterior cone conditions.

Recall that a cone $C$ is called {\it proper} if it does not contain any affine subspace of dimension one. Moreover, recall that for an open convex cone $\Gamma\subset\R^d$ its {\it dual cone} is defined as
\[\Gamma^\circ:=\{\xi\in\R^d;\,\forall\,y\in\Gamma:\,\langle y,\xi\rangle\geq 0\}.\]
For $\Gamma\neq\emptyset$ it is a closed proper convex cone in $\R^d$. On the other hand, every closed proper convex cone $C$ in $\R^d$ is the dual cone of a unique non-empty, open, convex cone which is given by
\[\Gamma:=\{y\in\R^d;\,\forall\xi\in C\backslash\{0\}:\,\langle y,\xi\rangle >0\}.\]
The proof can be done by the Hahn-Banach Theorem (cf.\ \cite[p.\ 257, vol.\ I]{Hoermander}). Therefore, we use the notation $\Gamma^\circ$ also for arbitrary closed convex proper cones. Moreover, from now on we assume all open convex cones $\Gamma$ to be non-empty.

As a first result we obtain from Proposition \ref{Hoermander reformulated} the next proposition which is an analogue result to \cite[Corollary 8.6.11, vol.\ I]{Hoermander}.

\begin{prop}\label{regularity in proper cones}
Let $\Gamma$ be an open proper convex cone in $\R^d$, $x_0\in\R^d$, and $P$ a non-constant polynomial. If for $\Omega:=x_0+\Gamma$ no hyperplane $H$ with $\sigma_P(H^\perp)=0$ intersects $\overline {\Omega}$ only in $x_0$, the following holds.

Each $u\in\mathscr{D}'(\Omega)$ with $P(D)u\in C^\infty(\Omega)$ which is $C^\infty$ outside a bounded subset of $\Omega$ already belongs to $C^\infty(\Omega)$.
\end{prop}

{\sc Proof.} Let $u\in\mathscr{D}'(\Omega)$ satisfy $P(D)u\in C^\infty(\Omega)$ and assume that $u$ is $C^\infty$ outside a bounded subset of $\Omega$. Since $\Gamma$ is a proper cone, there is a hyperplane $\pi$ intersecting $\Omega$ only in $x_0$. Let $H_\pi$ be a halfspace with boundary parallel to $\pi$ such that $\Omega_1:=\Omega\cap H_\pi\neq\emptyset$ is unbounded and $u|_{\Omega_1}\in C^\infty(\Omega_1)$. Denoting $\Omega_2:=\Omega$ we have convex sets $\Omega_1\subset \Omega_2$ and by the hypothesis, each hyperplane $H$ with $\sigma_P(H^\perp)=0$ and $H\cap \Omega_2\neq\emptyset$ already intersects $\Omega_1$. Proposition \ref{Hoermander reformulated} now gives $u\in C^\infty(\Omega)$.\hfill$\square$\\

The following proposition contains some elementary geometric results which will be used in the sequel.

\begin{prop}\label{geometric considerations}
\begin{itemize}
    \item[a)] If $C\subset\R^d$ is closed, convex, and unbounded, then for every $x\in C$ there is $\omega\in S^{d-1}$ such that $x+t\omega\in C$ for every $t\geq 0$.
    \item[b)] Let $\Gamma^\circ\neq\{0\}$ be a closed proper convex cone in $\R^d$ and $N\in S^{d-1}$. For $c\in\R$ let $H_c:=\{x\in \R^d;\langle x,N\rangle =c\}$. Then the following are equivalent.
    \begin{itemize}
        \item[i)] $H_0\cap \Gamma^\circ=\{0\}$.
        \item[ii)] $N\in\Gamma$ or $-N\in\Gamma$.
        \item[iii)] If $x\in\R^d$ and $H_c\cap(x+\Gamma^\circ)\neq\emptyset$ then $H_c\cap(x+\Gamma^\circ)$ is bounded.
        \item[iv)] If $x\in H_c$ then $H_c\cap(x+\Gamma^\circ)=\{x\}$.
    \end{itemize}
\end{itemize}
\end{prop}

{\sc Proof.} Part $a)$. Let $x\in C$. Replacing $C$ by $C-x$ we may assume without loss of generality that $x=0$. Let $(x_n)_{n\in\N}$ be a sequence in $C$ with $|x_n|\geq n$ for all $n\in\N$. Because $0\in C$ we have $1/|x_n|\, x_n\in C$ for every $n\in\N$. Passing to a subsequence if necessary, we can assume that $(1/|x_n|\,x_n)_{n\in\N}$ converges to $\omega\in S^{d-1}$. For every $t\geq 0$ we have $t/|x_n|<1$ for $n$ sufficiently large, hence $t/|x_n|\,x_n\in C$ for $0\in C$ and $C$ is convex. Since $C$ is closed it follows that $t\omega\in C$.

Part b). By use of a translation and an appropriate change of the value $c$, we can assume throughout the proof that $x=0$. Obviously, i) is then equivalent to iv).

To show that i) implies ii) let
\[H^+:=\{x;\,\langle x,N\rangle >0\}\mbox{ and }H^-:=\{x;\,\langle x,N\rangle <0\}.\]
If $H^+\cap\Gamma^\circ\neq\emptyset$ then $H^-\cap\Gamma^\circ=\emptyset$. Indeed, assume there are $x\neq y$ in $\Gamma^\circ$ such that $\langle x,N\rangle >0$ and $\langle y,N\rangle <0$. Convexity of $\Gamma^\circ$ and $H_0\cap\Gamma^\circ=\{0\}$ imply the existence of $\lambda\in (0,1)$ such that $\lambda x+(1-\lambda)y=0$, hence $-x=(1-\lambda)/\lambda\,y$. Since $\Gamma^\circ$ is a cone and $(1-\lambda)/\lambda>0$ it follows that $-x\in\Gamma^\circ$. Hence $\{0\}\neq \mbox{span}\{x\}\subset\Gamma^\circ$ contradicting that $\Gamma^\circ$ is proper.

Analogously one shows that $H^-\cap\Gamma^\circ\neq\emptyset$ implies $H^+\cap\Gamma^\circ=\emptyset$. Moreover, assuming $H^+\cap\Gamma^\circ=\emptyset$ as well as $H^-\cap\Gamma^\circ=\emptyset$ implies $\Gamma^\circ\subset H_0$. This yields $\Gamma^\circ=\{0\}$ because of $\Gamma^\circ\cap H=\{0\}$, contradicting $\Gamma^\circ\neq\{0\}$.

Without loss of generality we therefore may assume that $H^+\cap\Gamma^\circ\neq\emptyset$. From the above we obtain $\Gamma^\circ\subset\{x;\,\langle x,N\rangle\geq 0\}$. Since $H\cap\Gamma^\circ=\{0\}$ it follows that for all $x\in\Gamma^\circ\backslash\{0\}$ we have $\langle x,N\rangle >0$ which shows ii).

That ii) implies i) is trivial.

In order to show that iii) implies i) assume that $H_0\cap\Gamma^\circ\neq\{0\}$. Then, there is $\omega\in S^{d-1}$ such that $t\omega\in H_0\cap\Gamma^\circ$ for every $t\geq 0$. If $x\in H_c\cap \Gamma^\circ$ it follows that $x+t\omega\in H_c$. Moreover, because of $x\in \Gamma^\circ$ we have
\[\forall\,y\in\Gamma, t\geq 0:\,\langle y,x+t\omega\rangle=\langle y,x\rangle +t\langle y,\omega\rangle\geq 0,\]
hence $x+t\omega\in H_c\cap\Gamma^\circ$ for all $t\geq 0$ contradicting the boundedness of $H_c\cap\Gamma^\circ$.

To show that i) implies iii) assume that $H_c\cap\Gamma^\circ\neq\emptyset$ is unbounded. It follows from $a)$ that for $x\in H_c\cap\Gamma^\circ\backslash\{0\}$ there is $\omega\in S^{d-1}$ such that $x+t\omega\in H_c\cap\Gamma^\circ$ for all $t\geq 0$. Thus
\[c=\langle x,N\rangle=\langle x,N\rangle +t\langle\omega,N\rangle,\]
i.e.\ $\omega\in H_0$, and
\[\forall y\in\Gamma,t\geq 0:\,0\leq\langle y,x+t\omega\rangle.\]
Since $\Gamma$ is a cone, this implies
\[\forall y\in\Gamma,t\geq 0,\varepsilon>0:\,0\leq\langle\varepsilon y,x+t/\varepsilon\,\omega\rangle=\varepsilon\langle y,x\rangle+t\langle y,\omega\rangle.\]
The special case $t:=\langle y,x\rangle$ gives
\[\forall y\in\Gamma, \varepsilon>0:\,0\leq (\varepsilon +\langle y,\omega\rangle)\langle y,x\rangle.\]
Because $x\in\Gamma^\circ\backslash\{0\}$ we have $\langle y,x\rangle>0$ for every $y\in\Gamma$, so that the above inequality yields $\langle y,\omega\rangle\geq 0$ for all $y\in\Gamma$, thus $\omega\in\Gamma^\circ$. We conclude that $\omega\in H_0\cap\Gamma^\circ\cap S^{d-1}$ contradicting i).\hfill$\square$\\

We are now able to prove the main result of this section.

\begin{theo}\label{p-convexity by cones}
Let $\Omega$ be an open connected subset of $\R^d$ and $P\in\C[X_1,\ldots,X_d]$ a non-constant polynomial with principal part $P_m$.
\begin{itemize}
    \item[i)] $\Omega$ is $P$-convex for supports if for every $x\in\partial \Omega$ there is an open convex cone $\Gamma$ such that $(x+\Gamma^\circ)\cap\Omega=\emptyset$ and $P_m(y)\neq 0$ for all $y\in\Gamma$.
    \item[ii)] $\Omega$ is $P$-convex for singular supports if for every $x\in\partial \Omega$ there is an open convex cone $\Gamma$ such that $(x+\Gamma^\circ)\cap \Omega=\emptyset$ and $\sigma_P(y)\neq 0$ for all $y\in\Gamma$.
\end{itemize}
\end{theo}

{\sc Proof.} The proofs of both parts are very similar, so we give the proof of part ii) and only sketch the proof of i).

Let $u\in\mathscr{E}'(\Omega)$. We set $K:=\singsupp P(-D)u$ and $\delta:=dist(K, \Omega^c)$. If we show that $\dist(\singsupp u,\Omega^c)=\delta$ it follows from \cite[Theorem 10.7.3, vol.\ II]{Hoermander} that $\Omega$ is $P$-convex for singular supports. Since $\singsupp u\supset \singsupp P(-D)u$ we only have to show that $dist(\singsupp u, \Omega^c)\geq\delta$.

Let $x_0\in\partial \Omega$ and let $\Gamma$ be as in the hypothesis for $x_0\in\partial \Omega$. Then $(x_0+\Gamma^\circ)\cap \Omega=\emptyset$, thus $(x_0+y+\Gamma^\circ)\cap K=\emptyset$ for all $y\in\R^d$ with $|y|<\delta$. Therefore, for fixed $y$ with $|y|<\delta$, there is an open proper convex cone $\tilde{\Gamma}$ in $\R^d$ with $\tilde{\Gamma}\supset\Gamma^\circ\backslash\{0\}$ such that $(x_0+y+\tilde{\Gamma})\cap K=\emptyset$. Hence, $u\in\mathscr{E}'(\Omega)\subset\mathscr{D}'(x_0+y+\tilde{\Gamma})$ satisfies $P(-D)u\in C^\infty(x_0+y+\tilde{\Gamma})$.

We will show that $u\in C^\infty(x_0+y+\tilde{\Gamma})$ by applying Proposition \ref{regularity in proper cones}. Hence, let $H=\{v\in\R^d;\langle v,N\rangle=\alpha\}$ be a hyperplane with $\sigma_P(N)=0$. As $\overline{\tilde{\Gamma}}$ is a closed proper convex cone with non-empty interior, it is the dual cone of some open proper convex cone $\Gamma_1$. It follows from $\Gamma_1^\circ=\overline{\tilde{\Gamma}}\supset\Gamma^\circ$ that $\Gamma_1\subset\Gamma$. Because $\sigma_P(N)=0$ it follows from the hypothesis that $\{N,-N\}\cap\Gamma=\emptyset$, hence $\{N,-N\}\cap\Gamma_1=\emptyset$, so that by Proposition \ref{geometric considerations} b)  $H$ does not intersect $x_0+y+\overline{\tilde{\Gamma}}$ only in $x_0+y$. Since $u\in\mathscr{E}'(\Omega)$ we have that $\singsupp u$ is compact. Moreover $P(-D)u\in C^\infty(x_0+y+\tilde{\Gamma})$, so that $u\in C^\infty(x_0+y+\tilde{\Gamma})$ by Proposition \ref{regularity in proper cones}.

Since $x_0\in\partial \Omega$ and $y$ with $|y|<\delta$ were chosen arbitrarily, it follows that $dist(\singsupp u, \Omega^c)\geq\delta$, which proves ii).

In order to prove i), let $u\in\mathscr{E}'(\Omega)$, $K:=\supp P(-D)u$ and $\delta:=dist(K,\Omega^c)$. By \cite[Theorem 10.6.3, vol.\ II]{Hoermander} one has to show $dist(\supp u, \Omega^c)\geq\delta$ which is done as in the proof of ii) by using \cite[Corollary 8.6.11, vol.\ I]{Hoermander} instead of Proposition \ref{regularity in proper cones}.\hfill$\square$

\section{Proof of Theorem \ref{treves conjecture}}

Recall that for elliptic $P$ every open subset $\Omega\subset\R^d$ is $P$-convex for supports. In case of $d=2$ a complete characterization of $P$-convexity for supports is known. It is due to H\"ormander, see e.g.\ \cite[Theorem 10.8.3, vol.\ II]{Hoermander}.

\begin{theo}\label{non-elliptic p-convexitiy for supports in r^2}
If $P$ is non-elliptic then the following conditions on an open connected set $\Omega\subset \R^2$ are equivalent.
\begin{itemize}
    \item[i)] $\Omega$ is $P$-convex for supports.
    \item[ii)] The intersection of every characteristic hyperplane with $\Omega$ is convex.
    \item[iii)] For every $x_0\in\partial \Omega$ there is a closed proper convex cone $\Gamma^\circ\neq\{0\}$ with $(x_0+\Gamma^\circ)\cap \Omega=\emptyset$ and no characteristic hyperplane intersects $x_0+\Gamma^\circ$ only in $x_0$.
\end{itemize}
\end{theo}

In view of Proposition \ref{geometric considerations} the above condition iii) clearly is equivalent to the following condition.
\begin{itemize}
    \item[iii')] {\it For every $x_0\in\partial \Omega$ there is an open convex cone $\Gamma\neq\R^2$ with $(x_0+\Gamma^\circ)\cap \Omega=\emptyset$ and $P_m(y)\neq 0$ for all $y\in\Gamma$, where $P_m$ denotes the principal part of $P$.}
\end{itemize}

An analogous theorem to Theorem \ref{non-elliptic p-convexitiy for supports in r^2} for $P$-convexity for singular supports is the following. Recall that by Remark \ref{hypoelliptic remark} a polynomial $P$ is hypoelliptic if and only if $\sigma_P(H^\perp)\neq 0$ for every hyperplane $H$.

\begin{theo}\label{non-hypoelliptic p-convexitiy for singular supports in r^2}
If $P$ is non-hypoelliptic then the following conditions on an open connected set $\Omega\subset \R^2$ are equivalent.
\begin{itemize}
    \item[i)] $\Omega$ is $P$-convex for singular supports.
    \item[ii)] The intersection of $\Omega$ with every hyperplane $H$ satisfying $\sigma_P(H^\perp)=0$ is convex.
    \item[iii)] For every $x_0\in\partial \Omega$ there is an open convex cone $\Gamma\neq\R^2$ with $(x_0+\Gamma^\circ)\cap \Omega=\emptyset$ and $\sigma_P(y)\neq 0$ for all $y\in\Gamma$.
\end{itemize}
\end{theo}

The proof of the above theorem follows almost exactly the same lines as the proof of \cite[Theorem 10.8.3, vol.\ II]{Hoermander}.

Recall that a real valued function $f$ defined on a subset $M$ of $\R^d$ is said to {\it satisfy the minimum principle in the closed subset $F$ of $\R^d$} if for every compact subset $K\subset F\cap M$ it holds that $\inf_{x\in K}f(x)=\inf_{x\in\partial_F K}f(x)$, where $\partial_F K$ denotes the boundary of $K$ relative $F$. Moreover, we denote by
\[d_\Omega:\Omega\rightarrow\R,x\mapsto\dist(x,\Omega^c)\]
the so called boundary distance.\\

{\sc Proof of Theorem \ref{non-hypoelliptic p-convexitiy for singular supports in r^2}.} i)$\Rightarrow$ ii) It is enough to show that if $(\pm
1,0)\in \Omega$ and $\sigma_P((0,1))=0$ (i.e.\ parallels to the $x$-axis
are hyperplanes $H$ with $\sigma_P(H^\perp)=0$), then
$I=[-1,1]\times\{0\}\subset \Omega$. We join $(-1,0)$ and $(1,0)$ by
a polygon $\gamma$ in $\Omega$ without self-intersection, where we
can assume that $\gamma$ intersects the $x$-axis only at its end
points. For if this is not the case we can decompose $\gamma$ into
several polygons meeting the $x$-axis only at the end points and
treat them separately. Then $I$ and $\gamma$ are the boundary of a
connected and compact set $C$. We define
\[Y=\{y;\, (x,y)\in C\mbox{ for some }x\}\]
\[Y_0=\{y\in Y;\, (x,y)\in C\Rightarrow (x,y)\in \Omega\}.\]
$Y$ is a closed interval with non-empty interior and $Y_0$ is not empty since the end point of $Y$ which is different from $0$ belongs to $Y_0$. Since $\Omega$ is $P$-convex for singular supports it follows from \cite[Corollary 11.3.2]{Hoermander} that $d_\Omega$ satisfies the minimum principle in the hyperplane $\R\times\{y\}$ for arbitrary $y\in\R$. Therefore, if $y\in Y_0$ then from the definition of $Y_0$ $(x,y)\in C$ implies $(x,y)\in \Omega$ so that $\emptyset\neq C\cap (\R\times\{y\})\subset \Omega\cap (\R\times\{y\})$ is compact. Hence for $y\in Y_0$ and $x$ with $(x,y)\in C$ we have due to the minimum principle
\[d_\Omega(x,y)\geq d_\Omega(C\cap (\R\times\{y\}))=d_\Omega(\partial C\cap (\R\times\{y\}))\geq d_\Omega(\gamma\cap (\R\times\{y\}))\geq d_\Omega(\gamma).\]
Since $\gamma\subset \Omega$ we have that $d_\Omega(\gamma)>0$,
i.e.\ if $y\in Y_0$ then $(x,y)\in C$ implies that the distance form
$(x,y)$ to $\Omega^c$ is bounded below by the positive constant
$d_\Omega(\gamma)$. From this it follows that $Y_0$ is closed in
$Y$. Since $\Omega$ is open $Y_0$ is also open in the interval $Y$.
$Y_0$ being not empty now implies that $Y=Y_0$, hence $0\in Y=Y_0$,
so that $I=[-1,1]\times\{0\}\subset \Omega$.

ii)$\Rightarrow$ iii) If $x_0\in\partial \Omega$ and $H$ is a
hyperplane through $x_0$ with $\sigma_P(H^\perp)=0$ then one half
ray $H_1$ of $H$ bounded by $x_0$ is contained in $\Omega^c$ by ii).
If there is another hyperplane $I$ through $x_0$ with $\sigma_P(I^\perp)=0$ such that $H_1\cap
I=\{x_0\}$ then one of its half rays $I_1$ bounded by $x_0$ is
contained in $\Omega^c$ by ii) and since $\Omega$ is connected it
can be chosen so that the convex hull $\Gamma^\circ$ of $H_1$ and
$I_1$ is contained in $\Omega^c$ (and obviously is a proper convex
cone by $H_1\cap I=\{x_0\}$). If there is a hyperplane $K$ through
$x_0$ with $\sigma_P(K^\perp)=0$ and with
$K\cap\Gamma^\circ=\{x_0\}$ we continue extending $\Gamma^\circ$
until there is no hyperplane $L$ with $\sigma_P(L^\perp)=0$
intersecting $\Gamma^\circ$ only in $x_0$. Observe that by Lemma
\ref{singulars are characteristic in r2} and the remark following it
this procedure stops after a finite number of extensions so that the
resulting closed convex cone is indeed proper! From Proposition
\ref{geometric considerations} it follows that for no $y\in\Gamma$
we have $\sigma_P(y)=0$.

$iii)\Rightarrow i)$ This follows from Theorem \ref{p-convexity by cones} $b)$ which itself was very much inspired by the proof of the corresponding implication of \cite[Theorem 10.8.3, vol.\ II]{Hoermander}.\hfill$\square$\\

The proof of Theorem \ref{treves conjecture} is now obvious.\\

{\sc Proof of Theorem \ref{treves conjecture}.} Without loss of
generality we can assume that $P$ is not hypoelliptic, hence not
elliptic. Moreover, by passing to the different components of
$\Omega$ we can assume without loss of generality that $\Omega$ is
connected.

As $\Omega$ is supposed to be $P$-convex for supports it follows from Theorem \ref{non-elliptic p-convexitiy for supports in r^2} that for every $x\in\partial\Omega$ there is a non-empty, open convex cone $\Gamma$ different from $\R^2$ such that $(x+\Gamma^\circ)\cap\Omega=\emptyset$ and $P_m(y)\neq 0$ for all $y\in \Gamma$. From Lemma \ref{singulars are characteristic in r2} it follows that $\sigma_P(y)\neq 0$ for every $y\in\Gamma$ so that Theorem \ref{non-hypoelliptic p-convexitiy for singular supports in r^2} implies the $P$-convexity for singular supports of $\Omega$.\hfill$\square$\\

Combining Theorem \ref{p-convexity by cones} with Example \ref{wave in rd} gives an easy example that an analogous conclusion for $d>2$ is not true in general.

\begin{example}
\begin{rm}
Let $d>2$ and $P(x_1,\ldots,x_d)=x_1^2-x_2^2-\ldots -x_d^2$. Moreover, let $\Gamma:=\{x\in\R^d;\, x_d>(x_1^2+\ldots +x_{d-1}^2)^{1/2}\}$. Then $\Gamma$ is an open convex cone with $\Gamma^\circ=\overline{\Gamma}$. Set $\Omega:=\R^d\backslash\overline{\Gamma}$. Since $\{x\in\R^d;\, P_2(x)=0\}\cap\Gamma=\emptyset$ it follows easily from Theorem \ref{p-convexity by cones} i) that $\Omega$ is $P$-convex for supports.

We have seen in Example \ref{wave in rd} that $\sigma_P(e_d)=0$,
where $e_d=(0,\ldots,0,1)$ so that the hyperplane
$H=\{x\in\R^d;\,\langle x,e_d\rangle=-1\}$ satisfies
$\sigma_P(H^\perp)=\sigma_P(e_d)=0$. Taking
$K:=H\cap\{x\in\R^d;|x|\leq 2\}$ it is easily seen that $d_\Omega$
does not satisfy the minimum principle in the hyperplane $H$.
Therefore, by \cite[Corollary 11.3.2, vol.\ II]{Hoermander} $\Omega$
is not $P$-convex for singular supports.
\end{rm}
\end{example}

\noindent\textbf{Acknowledgement.} I want to thank L.\ Frerick as well as D.\ Vogt for inspiring discussions.

\begin{small}
{\sc Bergische Universit\"at Wuppertal,
FB Mathematik und Naturwissenschaften,
Gau{\ss}str.\ 20,
D-42097 Wuppertal, GERMANY}

{\it E-mail address: kalmes@math.uni-wuppertal.de}
\end{small}

\end{document}